\documentclass[12pt]{article}

\usepackage{amsmath}
\usepackage{theorem}
\usepackage{amssymb}
\usepackage{latexsym}
\usepackage{a4wide}

\newtheorem{thm}{Theorem}
\newtheorem{lem}[thm]{Lemma}
\newtheorem{prop}[thm]{Proposition}

\newtheorem{alg}[thm]{Algorithm}

{\theorembodyfont{\rmfamily} \newtheorem{exa}[thm]{Example}}
\newenvironment{rem}{\noindent{\bf Remark.}}{\newline}
\newenvironment{pf}{\noindent{\bf Proof.}}{\hbox{}\hfill $\Box$}

\newcommand{\Q}{\mathbb{Q}}

\newcommand{\C}{\mathbb{C}}

\newcommand{\barF}{\overline{F}}
\newcommand{\Z}{\mathbb{Z}}

\newcommand{\Lie}{\mathrm{\mathop{Lie}}}

\newcommand{\Der}{\mathrm{\mathop{Der}}}
\newcommand{\GL}{\mathrm{\mathop{GL}}}

\newcommand{\ad}{\mathrm{\mathop{ad}}}
\newcommand{\diag}{\mathrm{\mathop{diag}}}
\newcommand{\gl}{\mathfrak{\mathop{gl}}}
\newcommand{\mf}[1]{\mathfrak{#1}}
\newcommand{\g}{\mathfrak{g}}

\begin{document}

\title{Constructing algebraic groups from their Lie algebras}
\author{Willem A. de Graaf\\
Dipartimento di Matematica\\
Universit\`{a} di Trento\\
Italy\\
email: {\tt degraaf@science.unitn.it}}
\date{}
\maketitle

\begin{abstract}
A connected algebraic group in characteristic $0$ is uniquely determined
by its Lie algebra. In this paper an algorithm is given for constructing 
an algebraic group in characteristic $0$, given its Lie algebra. Using 
this an algorithm is presented for finding a maximal reductive subgroup and 
the unipotent radical of an algebraic group.
\end{abstract}

\section{Introduction}

Due to a wide variety of applications of algebraic groups there is considerable interest in doing
computations with such groups. For example, \cite{derksenetal}, \cite{gruseg} 
contain both applications of algebraic groups and computational methods. 
The most natural way to give an algebraic group is by a set of polynomial
equations. However, it is not easy to determine properties of the group
from these equations (for example, whether it is simple, solvable, or
nilpotent). A standard way to get around this is to study the Lie 
algebra of the group. This Lie algebra reflects many of the properties of
the group, and can therefore be used to get a picture of the structure of the 
group. However, by looking at the Lie algebra alone it is for instance not 
possible to construct subgroups, such as the unipotent radical, or a maximal
reductive subgroup. It is the aim of this paper to provide methods for 
constructing such subgroups by taking the inverse route. In other words, 
we describe algorithms for going back to the group from its Lie algebra.\par
For this a few restrictions are necessary. Firstly, in characteristic 
different from zero it is in general not true that to a given Lie algebra
corresponds a unique algebraic group. There can be more than one. So we restrict
to fields of characteristic zero. Second, since the Lie algebra of an algebraic
group is equal to the Lie algebra of the connected component of the identity,
we restrict to connected algebraic groups. We note that, also with these restrictions,
not every Lie algebra corresponds to an algebraic group. Lie algebras that do are
called algebraic. There are methods to decide whether or not a given Lie algebra 
is algebraic (cf. \cite{fiegra}). So we assume that we have given an algebraic Lie
algebra, and we want to construct the corresponding group.\par
Now we state the problem more precisely. Let $F$ be a field and let $\overline{F}$ be its 
algebraic closure. By $R(n,F)$ we denote the polynomial ring $F[x_{11},\ldots,x_{nn}]$
($n^2$ indeterminates). A subgroup $G\subset \GL(n,\overline{F})$ is said
to be algebraic if there is a set of polynomials $S\subset R(n,\barF)$
such that $G$ consists of the $g\in \GL(n,\overline{F})$ with
$f(g)=0$ for all $f\in S$. Let $G$ be an algebraic group, and let $I(G)\subset R(n,\barF)$ 
be the ideal of all polynomials vanishing on $G$. Then $G$ is said to be defined over 
$F$ if $I(G)$ is generated by polynomials in $R(n,F)$. \par
In this paper we take $F$ to be a field of characteristic $0$ ($\Q$ for example),
and consider algebraic groups that are defined over $F$. Since $F$ is perfect, any 
subgroup $\subset \GL(n,\overline{F})$ defined as the vanishing set of
a set of polynomials in $R(n,F)$ is automatically defined
over $F$ (cf. \cite{hum2}, \S 34 ). We call a set $S\subset R(n,F)$
such that $G$ consists of all $g$ where all elements of $S$ vanish, a set of
defining polynomials for $G$. (We do not assume that they generate $I(G)$.) \par
The Lie algebra $\g=\Lie(G)$ of an algebraic group $G\subset \GL(n,\barF)$ is a subalgebra
of $\gl(n,\barF)$ (the Lie algebra of all $n\times n$-matrices with coefficients
in $\barF$). Now if $G$ is defined over $F$, then $\g$ is a subalgebra of $\gl(n,F)$.
So the problem that we consider is the following. Given an algebraic Lie subalgebra $\g\subset
\gl(n,F)$, compute a set of defining polynomials (over $F$) for the unique connected
subgroup $G\subset \GL(n,\barF)$ such that $\g= \Lie(G)$. \par
The algorithm outlined in this paper consists of a few
subalgorithms. In Section \ref{sec:prodgrp} we describe
an algorithm that constructs the smallest algebraic group containing two
given algebraic subgroups $G_1,G_2\subset \GL(n,\barF)$. In Section \ref{sec:nilpotent}
an algorithm is given for constructing the algebraic group corresponding to
a Lie algebra that consists of nilpotent matrices. Let $X\in \g$; then 
by $G(X)$ we denote the smallest algebraic subgroup of $G$ such that its Lie algebra 
contains $X$. In Section \ref{sec:semsim} we give an algorithm for constructing $G(X)$
in case $X$ is semisimple. Now these algorithms together solve the problem.
Indeed, since $\g$ is algebraic, it is closed under Jordan decomposition.
So $\g$ has a basis consisting of elements that are either nilpotent or semisimple. 
So for every basis element $X$ we can construct $G(X)$. Let $X,Y$ be two basis
elements. Then with the algorithm of Section \ref{sec:prodgrp} we can construct the 
smallest algebraic group $H$ that contains
$G(X)$ and $G(Y)$. The Lie algebra of $H$ is generated by $X,Y$. Then we take
a third basis element $Z$, outside of this subalgebra, and form the smallest
algebraic group that contains $H$ and $G(Z)$. Continuing like this eventually
we find $G$. \par
In Section \ref{sec:reductive} we give an algorithm to decompose the Lie
algebra $\g$ into a direct sum of a reductive subalgebra, and an ideal
consisting of nilpotent elements. Together with the algorithms of the
previous sections this yields an algorithm for finding the analogous
decomposition of an algebraic group. This is one of the steps needed
in the algorithms of \cite{gruseg}. Finally, the last section describes
experiences with an implementation of the algorithm in the computer algebra
system {\sc Magma}.

\section{The algebraic group generated by two subgroups}\label{sec:prodgrp}

Let $G_1,G_2\subset \GL(n,\barF)$ be two algebraic groups defined over $F$.
Let $G\subset \GL(n,\barF)$ be the smallest algebraic group containing
both $G_1$ and $G_2$. Then also $G$ is defined over $F$ (this follows for
example from the algorithm given below). We consider the problem of finding
a set of defining polynomials in $R(n,F)$ for $G$. \par
First we note that the multiplication map $\cdot : \GL(n,\barF)\times \GL(n,\barF) \to
\GL(n,\barF)$ is a morphism of algebraic varieties. Using Gr\"obner basis techniques 
we can compute the Zariski closure of the image of this morphism (see, e.g., 
\cite{greupfis}, \S 1.8.3). This algorithm only uses
operations in the base field. Hence the closure of the image of a variety defined over $F$
is defined over $F$ as well. So we can compute defining polynomials in $R(n,F)$ for the Zariski closure 
$\overline{G_1G_2}$ of the set $G_1G_2 = \{ g_1g_2 \mid g_i\in G_i \}$. 
Now from \cite{chevii} (Chapter II, \S 7, Corollary 3),  we have the following theorem.

\begin{thm}[Chevalley]\label{thm1.1}
Let $G_1,\ldots,G_s$ be connected algebraic subgroups of $\GL( n,  \barF )$. Then the 
group $H$ generated by $G_1,\ldots,G_s$ is algebraic and connected. Moreover, there 
exists an $m>0$ such that every element of $H$ can be written as a product of $m$ elements,
each belonging to a $G_i$. 
\end{thm}

\begin{prop}\label{prop1.2}
Let $G_1,G_2$ be connected algebraic subgroups of $\GL(n,\barF)$. Let $H\subset \GL(n,\barF)$
be the group generated by $G_1,G_2$. Let $A\subset H$ and suppose that $1\in A$ and
that $\overline{AG_1}=A$ and $\overline{AG_2}=A$. Then $A=H$. 
\end{prop}

\begin{pf}
Let $a\in A$, and $g\in G_1$. Then $ag\in \overline{AG_1}$, and consequently $ag\in A$.
In the same way we see that $ag\in A$ if $g\in G_2$. 
Let $h\in H$. Then by Theorem \ref{thm1.1}, $h$ can be written as a product of elements 
from $G_1\cup G_2$. It follows that $ah\in A$ for all $a\in A$ and $h\in H$. So since
$1\in A$ we conclude that $A=H$.
\end{pf}

This leads to the following algorithm for computing $G$:

\begin{alg}
Input: algebraic subgroups $G_1,G_2\subset \GL(n,\barF)$.\\
Output: defining polynomials for the smallest algebraic subgroup 
of $\GL(n,\barF)$ containing $G_1,G_2$. \\
\begin{enumerate}
\item Let $G$ be the trivial subgroup of $\GL(n,\barF)$. 
\item Set $G' := G$.
\item Set $G' := \overline{G'G_1}$.
\item Set $G' := \overline{G'G_2}$.
\item If $G' \neq G$ then set $G:= G'$ and return to 2. Otherwise
return $G$.  
\end{enumerate}
\end{alg}

\begin{pf}
We show that the algorithm is correct and that it terminates.
The identity is always contained in $G'$. Secondly, if in a round of the
iteration $G'$ does not change then $\overline{G'G_1}=G' = \overline{G'G_2}$. So
by Proposition \ref{prop1.2}, $G'=H$, where $H$ is the group generated by $G_1,G_2$.\par
After $k$ rounds of the iteration $G_1G_2\cdots G_1G_2$ ($k$ factors $G_1G_2$) is
a subset of $G$. So by Theorem \ref{thm1.1} there is an $m>0$ such that after $m$
rounds of the iteration we have $G=H$. But then $G'$ will not change in the next round of 
the iteration, and the algorithm terminates. 
\end{pf}

\section{The nilpotent case}\label{sec:nilpotent}

Let $N\subset \gl(n,F)$ be a matrix Lie algebra consisting of nilpotent matrices. 
Then $N$ is an algebraic Lie algebra (\cite{cheviii}, Chapter V, Proposition 14).
By $G(N)\subset \GL(n,\barF)$ we denote the connected algebraic group having Lie algebra 
$N$. We consider the problem of finding defining polynomials for $G(N)$. \par
Let $N_1,\ldots,N_r$ be a basis of $N$. Consider the polynomial ring 
$P=F[T_1,\ldots ,T_r ,x_{ij} \mid 1\leq i,j\leq n]$.
Form the matrix $M=\exp( T_1N_1+\cdots +T_rN_r)$, and let $I$ be the
ideal of $P$ generated by $x_{ij}-M(i,j)$. By elimination techniques using Gr\"obner bases
we can compute the ideal $J= I \cap F[x_{ij}]$. 

\begin{prop}
Set $J'=\{ f\in F[x_{ij}] \mid f(a) = 0 \text{ for all } a\in G(N)\}$. Then 
$J=J'$.
\end{prop}

\begin{pf}
First of all we note that every element of $G(N)$ can uniquely be written as $\exp( \sum_i 
\alpha_i N_i)$, where $\alpha_i \in \barF$ (\cite{cheviii}, Chapter V, Proposition 14). 
Now let $f\in J$ and $a=\exp( \sum_i \alpha_i N_i ) \in G(N)$. We show that $f(a)=0$.
Consider the ring homomorphism $\varphi : F[T_1,\ldots,T_r,x_{ij}] \to F[x_{ij}]$,
defined by sending $T_i$ to $\alpha_i$, and leaving the $x_{ij}$ alone. Then
$\varphi(h) = h$ for all $h\in J$. Also $\varphi(h)(a)=0$ for all $h\in I$, as this
holds for the generators of $I$. Hence $f(a) = \varphi(f)(a)=0$. It follows 
that $J\subset J'$.\par
Now let $f\in J'$, then $f(\exp( T_1N_1+\cdots +T_rN_r))=0$ because 
$f(\sum_i \alpha_i N_i)=0$ where the $\alpha_i$ are arbitrary elements of $F$.
Since $x_{ij}=M(i,j)+(x_{ij}-M(i,j))$, every polynomial in $F[x_{ij}]$ can be
written as the sum of a polynomial in $T_1,\ldots,T_r$ and an element of $I$.
In particular, $f = p(T_1,\ldots,T_r) +h$ for a $h\in I$. Then
we substitute $x_{ij}\mapsto M(i,j)$. It follows that $p=0$ and $f\in I$.
Hence $f\in J$.
\end{pf}

\begin{rem}
One of the advantages of this approach is that it returns the vanishing ideal
of $G(N)$. In \cite{gruseg} a more direct method is given. However,
it does not have the same advantage. It works as follows. Let $p_i\in R(n,F)$ 
for $1\leq i\leq m$ be the
linear polynomials that define $N$ as a linear subspace of $\gl(n,F)$.
Let $X=(x_{ij})$ be the matrix consisting of the indeterminates $x_{ij}$.
Then $G(N)$ is defined by the polynomial equations
$$(X-1)^n=0 \text{ and } p_i( \log^*(X) ) \text{ for } 1\leq i\leq m,$$
where $\log^*(X) = -\sum_{i=1}^{n-1} i^{-1} (1-X)^i$. 
However, the resulting polynomials are all of degree $n-1$, or $n$, and it can
prove very difficult to find generators of the vanishing ideal from them. 
For larger $n$ it is also rather difficult to construct these polynomials.
We constructed the following example. Let $\mf{g}$ be the simple Lie algebra
of type $A_3$. This Lie algebra has a $10$-dimensional irreducible representation.
Let $x_1,\ldots,x_6$ denote the matrices of the positive root vectors of $\mf{g}$
in this representation, and let $N$ be the space spanned by them.
We tried to construct the above polynomials in {\sc Magma}. But the system crashed when 
it exceeded the limit of 2GB of memory. The method based on elimination used
136 seconds to construct defining polynomials (see Section \ref{sec:practice}). 
\end{rem}

\section{The semisimple case}\label{sec:semsim}

Let $X\in \gl(n,F)$ be a semisimple matrix. By $G(X)$ we denote the smallest
algebraic subgroup of $\GL(n,\barF)$ such that its Lie algebra contains $X$. 
In this section we consider the problem of obtaining defining polynomials for $G(X)$.
We will describe an algorithm based on the following theorem.

\begin{thm}[Chevalley]\label{thm:chev}
Let $\alpha_i,\ldots,\alpha_n\in \barF$ be the eigenvalues of $X$. 
Let $A\in \GL(n,\barF)$ be such that 
$AXA^{-1}=\diag(\alpha_1,\ldots,\alpha_n)$. Denote this last matrix by $Y$.
Set
$$\Lambda = \{(e_1,\ldots,e_n)\in \Z^n \mid \sum_{i=1}^n \alpha_ie_i = 0\}.$$
Then  
$$G(Y) = \{ \diag(c_1,\ldots,c_n) \mid c_i\in \barF \text{ and } \prod_i c_i^{e_i}=1 
\text{ for all $(e_1,\ldots,e_n)\in \Lambda$} \}.$$
Furthermore, $G(X) = A^{-1}G(Y)A$.
\end{thm}

The first statement is \cite{chevii}, Chapter II, Proposition 2. The second statement
is immediate.\par
The first step in our algorithm will be to construct a finite extension $F'\supset F$
containing the eigenvalues $\alpha_i$. Then by linear algebra we can construct a 
basis of the space
$$\Lambda_\Q = \{(e_1,\ldots,e_n)\in \Q^n \mid \sum_{i=1}^n \alpha_ie_i = 0\}.$$
From this we want to get a basis of $\Lambda$. First of all, by multiplying by suitable
scalars we get a basis $B$ of $\Lambda_\Q$ whose elements have integral coordinates.
By $B$ we also denote the $m\times n$-matrix having the elements of $B$
as rows. Let $S$ be the Smith normal form of $B$ (cf. \cite{sims}). This
means that we have unimodular matrices $P,Q$ such that $PBQ=S$, with 
$$S=\begin{pmatrix} d_1 &  &  & 0 & \cdots & 0\\
                     & \ddots &  &  &  \\
                     & & d_m & 0 & \cdots & 0\\
\end{pmatrix},$$
where the $d_i$ are positive integers. 

\begin{lem}\label{lem_grp_1}
$B$ is a basis of $\Lambda$ if and only if all $d_i$ are equal to $1$.
\end{lem}

\begin{pf}
Suppose that all $d_i=1$. Let $\Lambda'\subset \Z^n$ be the subgroup
generated by $B$. Then there is a surjective homomorphism 
$f: \Z^n \to \Z^{n-m}$ with kernel $\Lambda'$ (\cite{sims}, Proposition
3.3 of Chapter 8). Let $v\in \Lambda$, and write $v$ as a linear combination
of elements of $B$ with rational coefficients. After multiplying by a suitable
integer $s$ we see that $sv\in \Lambda'$. But then $0=f(sv)=sf(v)$. Hence 
$f(v)=0$ and $v\in \Lambda'$. We conclude that $\Lambda=\Lambda'$.\par
Now suppose that $B$ is a basis of $\Lambda$. We have a surjective 
homomorphism $f : \Z^n \to \Z/d_1\Z\oplus \cdots \oplus\Z/d_m\Z  \oplus \Z^{n-m}$,
with kernel $\Lambda$. Suppose that $i$ is such that $d_i>1$. Then there is a
$v\in \Z^n$ such that $f(v)\neq 0$ but $f(d_iv)=0$. In other words, $d_iv\in \Lambda$.
But then $v\in \Lambda$, and we have a contradiction.  
\end{pf}

Note that $B=P^{-1}SQ^{-1}$. This means that $d_i$ times the $i$-th row of 
$Q^{-1}$ lies in the span of $B$, and hence in $\Lambda$. Therefore, the
$i$-th row of $Q^{-1}$ itself lies in $\Lambda$, for $1\leq i\leq m$. Set
$$S'=\begin{pmatrix} 1 &  &  & 0 & \cdots & 0\\
                     & \ddots &  &  &  \\
                     & & 1 & 0 & \cdots & 0\\
\end{pmatrix},$$
and $A=P^{-1}S'Q^{-1}$. Then the rows of $A$ belong to $\Lambda$, and form
a basis of $\Lambda_\Q$. Moreover, the Smith normal form of $A$ is $S'$. 
Hence by Lemma \ref{lem_grp_1}, the rows of $A$ form a basis of $\Lambda$.
We remark that there are efficient algorithms for computing the matrices
$P,Q$ (see \cite{sims}). Hence we can efficiently compute a basis of $\Lambda$,
given a basis of $\Lambda_\Q$.\par
Let $A_I(X)$ be the associative algebra with one generated by $X$. Then 
$A_I(X)$ is spanned by $I,X,X^2,\ldots,X^t$, where $t+1$ is the degree of the
minimal polynomial of $X$. By \cite{chevii}, \S 13, Theorem 10, $G(X)\subset A_I(X)$.
Let $\alpha_i,Y$ be as in Theorem \ref{thm:chev}. 
We note that the minimal polynomial of a semisimple matrix is the square free
part of its characteristic polynomial. In particular, it does not depend on the 
base field, but only on the coefficients of the matrix. So $A_I(Y)$ is spanned
by $I,Y,\ldots,Y^t$. \par
Let $y=\sum_{i=0}^t \delta_i Y^i\in A_I(Y)$. Then by $y(k,k)$ we denote the
entry on position $(k,k)$. By Theorem \ref{thm:chev} we have that $y\in G(Y)$
if and only if $\prod_k y(k,k)^{e_k} =1 $ for $(e_1,\ldots,e_n)$
in a basis of $\Lambda$. Now set $e_k'=e_k$ if $e_k\geq 0$,  and $e_k'=0$
otherwise. Also $e_k''=e_k'-e_k$. Then $\prod_k y(k,k)^{e_k} =1 $ if and only if 
$\prod_k y(k,k)^{e_k'} = \prod_k y(k,k)^{e_k''}$. In this we substitute
$y(k,k) = \sum_i \delta_i \alpha_k^i$. This yields a polynomial equation for the 
$\delta_i$ with coefficients in $F'$. By writing the coefficients
as linear combinations of a basis of $F'$ over $F$, we get polynomials 
$p_1,\ldots,p_s$, $p_k\in F[T_0,\ldots,T_t]$, with the property that
$\prod_k y(k,k)^{e_k} =1 $ if and only if $p_i(\delta_0,\ldots,\delta_t)=0$
for $1\leq i\leq s$. Let $h_1,\ldots,h_m$ be the totality of these polynomials 
that we get when we let $(e_1,\ldots,e_n)$ run through a basis of $\Lambda$.
Then $y\in G(Y)$ if and only if $h_i(\delta_0,\ldots,\delta_t)=0$
for $1\leq i\leq m$. \par
Note that $y\in G(Y)$ if and only if $A^{-1}yA\in G(X)$ (where $A$ is as 
in Theorem \ref{thm:chev}). But $A^{-1}yA = \sum_i \delta_i X^i$. We conclude that
$\sum_i \delta_i X^i \in G(X)$ if and only if $h_i(\delta_0,\ldots,\delta_t)=0$
for $1\leq i\leq m$. \par
We summarise these findings in the following algorithm:

\begin{alg}
Input: a semisimple matrix $X\in \gl(n,F)$.\\
Output: defining polynomials for $G(X)$. \\
\begin{enumerate}
\item Construct a finite extension $F'\supset F$ containing the eigenvalues
of $X$.
\item Compute a basis of $\Lambda$.
\item Construct polynomials $h_1,\ldots,h_m\in F[T_0,\ldots,T_t]$ with the
property that $\sum_i \delta_i X^i \in G(X)$ if and only if $h_i(\delta_0,\ldots,\delta_t)=0$
for $1\leq i\leq m$.
\item Set $M = \sum_{i=0}^t T_i X^i$ and let $\alpha_{ij}^k\in F$ be such that 
$T_k = \sum_{ij} \alpha_{ij}^k M(i,j)$. 
\item Let $x_{ij}$ for $1\leq i,j\leq n$ be indeterminates, and consider the
substitution $T_k \mapsto \sum_{ij} \alpha_{ij}^k x_{ij}$. Let $\varphi : 
F[T_0,\ldots,T_t] \to F[x_{ij}]$ be the corresponding ring homomorphism.
\item Let $g_1,\ldots,g_r\in F[x_{ij}]$ be linear polynomials with the 
property that $a\in \gl(n,F)$ lies in $A_I(X)$ if and only if $g_i(a)=0$ for
$1\leq i\leq r$.
\item Return $\{g_1,\ldots,g_r\}\cup \{ \varphi(h_1),\ldots,\varphi(h_m)\}$.
\end{enumerate}
\end{alg}

\begin{pf}
Steps 1, 2, and 3 have already been commented on. Since the $X^i$ for $0\leq i\leq t$
are linearly independent, the coefficients $\alpha_{ij}^k$ in Step 4 exist (but
they are not necessarily unique). Let $a=(a_{ij})\in \GL(n,F)$ be such that
$g_i(a)=0$ for all $i$. Then $a\in A_I(X)$ so $a= \sum_i \delta_i X^i$.
By substituting $T_k\mapsto \delta_k$ in the equation $T_k = \sum_{ij} \alpha_{ij}^k M(i,j)$
we see that $\delta_k = \sum_{i,j} \alpha_{ij}^k a_{ij}$. It follows that 
$a\in G(X)$ if and only if 
$$h_r( \sum_{ij}\alpha_{ij}^0 a_{ij}, \ldots, \sum_{ij}\alpha_{ij}^t a_{ij} ) = 0$$
for $1\leq r \leq m$. But this is equivalent to $\varphi(h_r)(a)=0$ for $1\leq r\leq m$.
We conclude that this algorithm returns defining polynomials for $G(X)$. 
\end{pf}

\begin{exa}
Let 
$$X= \begin{pmatrix} 0 & 1 \\ -1 & 0 \end{pmatrix}.$$
Then $A_I(X)$ is spanned by $I,X$. The eigenvalues of $X$ are $i,-i$, and
hence $\Lambda$ is spanned by $(1,1)$. Let $Y=\diag(i,-i)$ and $y=\delta_0 I
+\delta_1Y$, then $y(1,1) = \delta_0 +\delta_1 i$, $y(2,2)=\delta_0-\delta_1 i$.
Since $\dim \Lambda=1$, there is only one equation that needs to be satisfied, 
i.e., $y(1,1)y(2,2)=1$, or $\delta_0^2+\delta_1^2=1$.
Now 
$$M = T_0 I +T_1 X = \begin{pmatrix} T_0 & T_1 \\ 
-T_1 & T_0 \end{pmatrix}.$$
So we get the substitution $T_0 = x_{11}$ and $T_1=x_{12}$. 
The equation above then transforms to $x_{11}^2+x_{12}^2-1=0$.
Furthermore, the linear equations of Step 6 of the algorithm are $x_{21}+x_{12}=0$, 
$x_{11}-x_{22}=0$. 
\end{exa}


\section{An application}\label{sec:reductive}

Let $G\subset \GL(n,\barF)$ be an algebraic group, defined over $F$. 
Then there is a unique maximal unipotent subgroup,
$u(G)$, called the unipotent radical of $G$. Furthermore, there is a reductive
subgroup $H\subset G$ with $G=H\ltimes u(G)$. In \cite{gruseg} algorithms are given
for finding defining polynomials for $H$ and $u(G)$. However, these algorithms
do not seem to be very practical. In this section we describe a method based on
the algorithms of the previous sections.  \par
Set $\mathfrak{g}=\Lie(G)$. Let $\mathfrak{s}$ denote the solvable radical of 
$\mathfrak{g}$, and $\mathfrak{l}$ a Levi subalgebra of $\mathfrak{g}$, and $\mathfrak{n}$
the largest ideal of $\mathfrak{g}$ consisting of nilpotent elements. 
Then by \cite{cheviii}, Chapter V, \S 4, Proposition 5 (see also \cite{mostow}),
$\mathfrak{s}$ has a commutative subalgebra $\mathfrak{d}$ consisting of semisimple
elements, with the following properties
\begin{enumerate}
\item $\mathfrak{s}=\mathfrak{d}+\mathfrak{n}$ (semidirect sum),
\item $[\mathfrak{l},\mathfrak{d}]=0$,
\item $\mathfrak{n}$ is the set consisting of all nilpotent elements of $\mathfrak{s}$.
\end{enumerate}

Now $H$ and $u(G)$ are the algebraic groups corresponding to $\mathfrak{l}+\mathfrak{d}$ and
$\mathfrak{n}$ respectively. There are algorithms to compute $\mathfrak{l}$ and $\mathfrak{s}$
(cf. \cite{gra6}). In the remainder of this section we describe how to
compute $\mathfrak{d}$ and $\mathfrak{n}$. Then by using the algorithms in the first
part of this paper, we can compute defining polynomials for $H$ and $u(G)$. 

\begin{lem}\label{lem4.1}
Set $\mf{h}= C_\mf{s}( \mf{d} )$, the centralizer of $\mf{d}$ in $\mf{s}$. Then $\mf{h}$ is
a Cartan subalgebra of $\mf{s}$, and $\mf{h} = \mf{d} + C_{\mf{n}}(\mf{d})$.
\end{lem}

\begin{pf}
Let $x\in C_\mf{s}( \mf{d} )$, then we can write $x = y+z$, where $y\in \mf{d}$ and 
$z\in \mf{n}$. So for any $u\in \mf{d}$ we get $0=[u,x] = [u,y]+[u,z]=[u,z]$. It follows
that $z \in C_\mf{s}( \mf{d} )$. Hence $\mf{h} = \mf{d} + C_{\mf{n}}(\mf{d})$. Since
$\mf{n}$ is nilpotent and $[\mf{d},C_{\mf{n}}(\mf{d})]=0$
it also follows that $\mf{h}$ is nilpotent.\par
Now let $x\in \mf{d}$. Then $x$ is a semisimple linear transformation. This implies that 
also $\ad x :  \mf{s}\to \mf{s}$ is semisimple. Hence $C_\mf{s}( \mf{d} ) = \mf{s}_0( \ad \mf{d} )$,
where
$$\mf{s}_0( \ad \mf{d} ) = \{ y\in \mf{s} \mid (\ad x)^{\dim \mf{s}} (y)=0 
\text{ for all } x\in \mf{d} \}.$$
By \cite{win}, Theorem 4.4.4.8 every Cartan subalgebra of 
$\mf{s}_0( \ad \mf{d} )$ is a Cartan subalgebra of $\mf{s}$. But as seen above
$\mf{h}=\mf{s}_0( \ad \mf{d} )$ is nilpotent. So it is its own Cartan subalgebra. 
\end{pf}

\begin{lem}\label{lem4.2}
Let $\mf{h}'$ be a Cartan subalgebra of $\mf{s}$. Then there is an abelian subalgebra
$\mf{d}'\subset \mf{s}$ consisting of semisimple elements, such that $\mf{s} = \mf{d}'
+\mf{n}$ (semidirect sum), and $\mf{h}' = C_\mf{s}(\mf{d}') = \mf{d}' + C_\mf{n}(\mf{d}')$. 
\end{lem}

\begin{pf}
Let $D$ be the group of automorphisms of $\mf{s}$ generated by $\exp( \ad x)$ for 
$x\in [\mf{s},\mf{s}]$. Set $\mathfrak{h}= C_{\mf{s}}(\mf{d})$, which is a Cartan
subalgebra of $\mf{s}$ by Lemma \ref{lem4.1}. By \cite{cheviii}, Chapter VI, Proposition 19, 
there is a $\sigma\in D$ such that $\sigma(\mf{h}) =
\mf{h}'$. Set $\mf{d}' = \sigma( \mf{d} )$. In order to prove that $\mf{d}'$ consists
of semisimple elements we may assume that $\sigma = \exp( \ad x)$. Then $\sigma(u)
= (\exp x) u (\exp x)^{-1}$ (by \cite{jac}, Chapter IX, (38), note that 
$[\mf{s},\mf{s}]\subset \mf{n}$, so $\exp(x)$
is a well defined endomorphism). So $\mf{d}'$ is obtained from $\mf{d}$ by conjugation
by a fixed element, and hence it consists of semisimple elements. The other properties
follow from the analogous properties of $\mf{h}$ and the fact that $\sigma$ is an automorphism. 
\end{pf}

Let $\mf{r}$ be a finite-dimensional Lie algebra, and
let $\mf{a}\subset \mf{r}$ be a nilpotent subalgebra. Then $\mf{r}$ has a Fitting
decomposition with respect to the adjoint action of $\mf{a}$ (cf. \cite{jac}, \cite{gra6}).
This decomposition is written $\mf{r} = \mf{r}_0 (\mf{a}) \oplus \mf{r}_1(\mf{a})$. 
These are called the Fitting 0-component and Fitting 1-component respectively. We have
$\mf{r}_1(\mf{a}) = \cap_{i>0} [ \mf{a}^i, \mf{r} ]$, where by $[\mf{a}^i,\mf{r}]$ we denote
the space $[\mf{a},[\mf{a},\cdots,[\mf{a},\mf{r}]\cdots]]$ ($i$ factors $\mf{a}$). 
Based on this there is a 
straightforward algorithm  for computing a basis of  $\mf{r}_1(\mf{a})$ (see \cite{gra6}).

\begin{lem}
Let $\mf{h}'$ and $\mf{d}'$ be as in the previous lemma. Let $a_1,\ldots,a_r$ be a basis
of $\mf{h}'$ and let $a_i=s_i+n_i$ be the Jordan decomposition of $a_i$. Then $\mf{d}'$ is
spanned by $s_1,\ldots,s_r$. Let $\mf{s}_1(\mf{h}')$ be the Fitting 1-component of $\mf{s}$
with respect to the adjoint action of $\mf{h}'$. Then $\mf{n}$ is spanned by the $n_i$ along
with $\mf{s}_1(\mf{h}')$.
\end{lem}

\begin{pf}
Let $h\in \mf{h}'$; then by Lemma \ref{lem4.2}, there are $s\in \mf{d}'$ and $n\in 
C_\mf{n}(\mf{d}')$ with $h=s+n$. But then $s$ is semisimple, $n$ is nilpotent and
$[s,n]=0$. It follows that $h=s+n$ is the Jordan decomposition of $h$. In particular,
all $s_i$ lie in $\mf{d}'$ and all $n_i$ lie in $C_\mf{n}(\mf{d}')$. Also 
$h=\sum_i \alpha_i a_i = (\sum_i \alpha_i s_i) + (\sum_i \alpha_i n_i)$.
This implies that $\sum_i \alpha_i s_i = s$ and $\sum_i \alpha_i n_i = n$. 
Therefore, the $s_i$ span $\mf{d}'$, and the $n_i$ span $C_\mf{n}(\mf{d}')$.
Note that $\mf{s}_0(\mf{h}') = \mf{h}'$ as $\mf{h}'$ is a Cartan subalgebra of $\mf{s}$
(cf. \cite{jac}, Chapter III, Proposition 1).  
Also $\mf{s}_1(\mf{h}')\subset [\mf{s},\mf{s}] \subset \mf{n}$. From this we
get the last statement.
\end{pf}

We note that there are efficient algorithms for computing a Cartan subalgebra of
a Lie algebra of characteristic $0$ (cf. \cite{gra6}). So  the previous lemma 
yields a straightforward algorithm for computing bases of $\mf{d'}$ and $\mf{n}$.

\section{Practical experiences}\label{sec:practice}

In this section we report on practical experiences with an implementation
of the algorithms in the computer algebra system {\sc Magma}. The computations were
done on a 2GHz processor, with 2GB of memory. \par
First we consider the algorithm from Section \ref{sec:nilpotent}.
This algorithm takes as input a set of nilpotent matrices that span a 
Lie algebra, and returns defining polynomials for the corresponding
algebraic group. In order to generate input we have constructed irreducible
representations of several simple Lie algebras, and taken the matrices
of the positive root vectors. The running times of the algorithm are listed
in Table \ref{tab1}. From the table we see that the algorithm performs 
quite well for small matrices. However, the running times increase sharply
when the size of the matrices increases. This is due to the Gr\"obner basis 
calculation, used for the elimination of variables.

\begin{table}
\begin{center}
\begin{tabular}{|l|r|r|r|}
\hline
type & dim rep. & $\dim N$ & time $G(N)$ \\
\hline
$B_2$ & 5 & 4 & 0.08 \\
$B_2$ & 10 & 4 & 43 \\
$G_2$ & 7 & 6 & 3.45 \\
$A_3$ & 10 & 6 & 136 \\
$C_3$ & 14 & 9 & $\infty$\\
$D_4$ & 8 & 12 & $\infty$\\
\hline
\end{tabular}
\end{center}
\caption{Running times for the algorithm of Section \ref{sec:nilpotent}.
The first column lists the type of the simple Lie algebra, the second
column has the dimension of the irreducible module, the third column the dimension
of the Lie algebra spanned by the positive root vectors. The last column displays 
the time (in seconds) for constructing defining polynomials for the group.
An $\infty$ means that the algorithm did not complete with 2GB of memory.}
\label{tab1}
\end{table}

Next we consider the algorithm of Section \ref{sec:semsim}. Here the
input is a single semisimple matrix. As input we have used the companion
matrices of several square free polynomials (except the first one these are taken
from the database in \cite{klumal}).
The running times are displayed in Table \ref{tab2}. As was to be expected the
algorithm spends a lot more time if the degree of the splitting field gets larger.

\begin{table}
\begin{center}
\begin{tabular}{|l|r|r|}
\hline
polynomial & dimension splitting field & time $G(X)$ \\
\hline
$x^8-x^7+x^6+2x^5-3x^4+4x^3+2$ & 360 & 83 \\
$x^6-2x^4+x^2-2x-1$ & 360 & 293 \\
$x^8+12x^6+50x^4+83x^2+43$ & 384 & 89 \\
$x^8-2x^6+7x^4-8x^2-4x+7$ & 576 & 1463 \\
\hline
\end{tabular}
\end{center}
\caption{Running times for the algorithm of Section \ref{sec:semsim}.
The first column lists the polynomial, and the second column the degree of
its splitting field. The last column displays the time (in seconds) for constructing 
defining polynomials for the group corresponding to the companion matrix of the
polynomial.}
\label{tab2}
\end{table}

Finally we consider the algorithm of Section \ref{sec:prodgrp}. For this
algorithm we have constructed several rather different inputs, which
we describe separately.\par
Let $A$ be an algebra, and $\Der(A)$ the Lie algebra of its derivations.
By \cite{chevii}, Chapter II, \S 14, Theorem 16, we know that $\Der(A)$ is
the Lie algebra of the automorphism group of $A$. Here we consider the
algebra $A$ of quaternions over $\Q$. This algebra has basis $1,i,j,k$,
where $1$ is the identity. The other product relations are $ij=k$, $ji=-k$,
$ik=-j$, $ki=j$, $jk=i$, $kj=-i$. Then $\Der(A)$ is spanned by 
$$x_1=\begin{pmatrix}
0 &  0 &   0 &   0 \\
0 &  0 &   1 &   0 \\
0 & -1 &   0 &   0 \\
0 &  0 &   0 &   0 
\end{pmatrix},
x_2=\begin{pmatrix}
0 &  0 &   0 &   0 \\
0 &  0 &   0 &   1 \\
0 &  0 &   0 &   0 \\
0 & -1 &   0 &   0 
\end{pmatrix},
x_3=\begin{pmatrix}
0 &  0 &   0 &   0 \\
0 &  0 &   0 &   0 \\
0 &  0 &   0 &   1 \\
0 &  0 &  -1 &   0 
\end{pmatrix},$$
(where we use the column convention to represent a linear map with respect to
the basis $1,i,j,k$ of $A$). Let $H_i$ be the smallest algebraic group containing
$x_i$. The $x_i$ are semisimple, so equations for the $H_i$ can be computed 
using the algorithm of Section \ref{sec:semsim}. Let $G$ be the smallest algebraic
subgroup of $\GL(4,\C)$ containing $H_1$, $H_2$. The Lie algebra of $G$ contains
$x_1,x_2$. But then it contains also $x_3$. It follows that its Lie algebra is
equal to $\Der(A)$. Hence $G$ is the connected component of the identity of the automorphism group
of $A$. The computation of equations for $G$ (given those of $H_1$ and $H_2$) cost 8.3 seconds. \par
The next example has the Lie algebra consisting of the matrices
$$\begin{pmatrix}
* &  * &   * &   * \\
* &  * &   * &   * \\
0 &  0 &   * &   * \\
0 &  0 &   * &   * 
\end{pmatrix}.$$

Let $H_1$ denote the group corresponding to the subalgebra consisting of strictly
upper triangular matrices. By $H_{i,j}$ we denote the group corresponding to the 
subalgebra spanned by $e_{i,j}$ (where $e_{ij}$ denotes the matrix with $1$ on
position $(i,j)$ and zeros elsewhere). Let $G_1$ denote the smallest algebraic
group containing both $H_1$ and $H_{2,1}$. The construction of defining polynomials
for $G_1$ cost 0.99 seconds. Let $G_2$ denote the smallest algebraic group
containing $G_1$ and $H_{4,3}$. The computation of defining polynomials for $G_2$
cost 2.4 seconds. Let $G_3$ be the smallest algebraic group containing $G_2$
and $H_{1,1}$. The computation of defining polynomials for $G_3$ cost 5.4 seconds.
Let $G=G_4$ denote the smallest algebraic group containing $G_3$ and $H_{3,3}$.
The computation of defining polynomials for $G_4$ cost 15 seconds. We note that
$G$ is equal to the algebraic group corresponding to the whole Lie algebra.\par
Let $L$ be the Lie algebra of type $B_2$ in its $4$-dimensional representation.
Let $H_1$ and $H_2$ denote the groups corresponding to the subalgebra spanned by the positive
and negative root vectors respectively. Let $G$ be the smallest algebraic group
containing both $H_1$ and $H_2$. The computation of defining polynomials for $G$
cost 173 seconds.\par
Finally, let $L$ be the simple Lie algebra of type $A_3$ in its $4$-dimensional
representation. Let $H_1,H_2$ and $G$ be as in the previous example. The algorithm
did not manage to compute defining polynomials for $G$ within 2GB.\par
Concluding we can say that the algorithms do work for small examples. However,
due to the Gr\"obner basis computations, needing a lot of memory and running time,
the algorithms are not yet practical for larger examples.

\def\cprime{$'$} \def\cprime{$'$} \def\cprime{$'$} \def\cprime{$'$}

\end{document}